\documentclass[review]{elsarticle}
\pdfoutput=1

\usepackage{lineno,hyperref}
\usepackage{amsmath,amssymb}
\usepackage{xcolor}

\journal{Nonlinear Analysis}

\begin{document}

\begin{frontmatter}

\title{Inverse problems for a fractional conductivity equation}

\author{Giovanni Covi}
\address{Department of Mathematics, University of Jyv\"askyl\"a, Finland}
\address{\emph{\texttt{giovanni.g.covi@jyu.fi}}}

\begin{abstract}
This paper shows global uniqueness in two inverse problems for a fractional conductivity equation: an unknown conductivity in a bounded domain is uniquely determined by measurements of solutions taken in arbitrary open, possibly disjoint subsets of the exterior. {\color{black} Both the cases of infinitely many measurements and a single measurement are addressed.} The results are based on a reduction from the fractional conductivity equation to the fractional Schr\"odinger equation, and as such represent extensions of previous works. Moreover, a simple application is shown in which the fractional conductivity equation is put into relation with a long jump random walk with weights.
\end{abstract}

\begin{keyword}
Fractional conductivity equation \sep non-local operators \sep Inverse problems \sep Calder\'on problem 
\MSC[2010] 35R11, 35R30
\end{keyword}

\end{frontmatter}

\section{Introduction} \label{intro}
This paper introduces and studies a fractional conductivity equation, and establishes uniqueness and reconstruction results for related inverse problems. The main point of interest is a fractional version of the standard Calder\'on problem \cite{Ca80}, which requires to find the electrical conductivity of a medium from voltage and current measurements on its boundary.
\vspace{5mm}

Let $\Omega\subset\mathbb R^n$ be a bounded open set with a regular enough boundary (e.g., let $\partial\Omega$ be Lipschitz), representing a medium whose electrical properties must be studied. The Dirichlet problem for the conductivity equation asks to find a function $u$ satisfying 
$$\left\{\begin{array}{lr}
        \nabla\cdot(\gamma\nabla u) =0 & \text{in } \Omega\\
        u=f & \text{on } \partial\Omega
        \end{array}\right. \;, $$ \vspace{1mm}

\noindent where $f$ is some prescribed boundary value and $\gamma$ is the electrical conductivity of the medium. The boundary measurements are given by the Dirichlet-to-Neumann (or DN) map $\Lambda_\gamma : H^{1/2}(\partial\Omega)\rightarrow H^{-1/2}(\partial\Omega)\;,$ which is defined weakly using the bilinear form of the equation. The inverse problem consists in finding the function $\gamma$ in $\Omega$ from the knowledge of $\Lambda_\gamma$.

\vspace{5mm}

{\color{black} The classical Calder\'on problem we stated above has general mathematical interest, as it serves as a model case for the study of inverse problems for elliptic equations, and is of course useful in the applied fields of medical, seismic and industrial imaging. The survey \cite{Uh14} provides many more details on this topic. The main physical motivation, and actually Calder\'on's original one, comes from electrical mineral prospecting. In this application, the electrical properties of a patch of soil are measured by an array of electrodes distributed on the ground, with the goal of determining whether any economically interesting mineral source is present underneath. }

\vspace{1mm}

\noindent {\color{black} On the other hand, fractional mathematical models are nowadays widely used in many fields of science. It is known for example that they arise in the study of turbulent fluids such as the atmosphere. They also appear in probability theory as generators of certain Levy processes, and because of this they are used in mathematical finance. For the many modern applications of fractional models, check e.g. \cite{BV18}. }

\vspace{1mm}

\noindent {\color{black} It is therefore very promising to study a fractional extension of the Calder\'on problem, in view of its many potential applications. This is the model we introduce below. }

\vspace{5mm}

Fix $s\in(0,1)$ and consider the new operators $(\nabla\cdot)^s$ and $\nabla^s$, which in this paper are called \emph{fractional divergence} and \emph{fractional gradient}. Their rigorous definitions will be given in Section \ref{prel} following \cite{DGLZ2013}, but for now they can be thought of as non-local counterparts of the standard divergence and gradient. They are "nonlocal" because they do not preserve supports, in the sense that $\nabla^su|_\Omega$ can only be computed knowing $u$ over all of $\mathbb R^n$. Later on we will show that, just as in the local case, the combination of these operators gives the fractional Laplacian, that is $(-\Delta)^s u = (\nabla\cdot)^s\nabla^s u$. 

\paragraph{Remark} It is worth noticing at this point that our choice for the names of the non-local operators, which has been guided by the similarity with the local case, is not universal. In \cite{DGLZ2013}, for example, our fractional gradient is called \emph{adjoint of the fractional divergence}, while the name \emph{fractional gradient} is assigned to a completely different operator which does not play any role in this paper. 
\vspace{5mm}

We set up the Dirichlet problem for the fractional conductivity equation as 

$$\left\{\begin{array}{lr}
        (\nabla\cdot)^s(\Theta\cdot\nabla^s u) =0 & \text{in } \Omega\\
        u=f & \text{in } \Omega_e
        \end{array}\right. \;,
				$$
\vspace{2mm}

\noindent where $\Theta$ is an interaction matrix depending on $\gamma$. Because of the non-local nature of the operators, the exterior value is given over all of $\Omega_e=\mathbb R^n\setminus \overline\Omega$. In section \ref{proof} it will be shown that the bilinear form associated to the conductivity equation is positive definite; this assures that 0 is not an eigenvalue of $(\nabla\cdot)^s(\gamma\nabla^s)$, and therefore the problem above is well-posed. Consequently, the DN map $\Lambda^s_\gamma : H^s(\Omega_e) \rightarrow (H^{s}(\Omega_e))^*$ can be defined in a weak sense starting from the bilinear form of the equation. The inverse problem asks to recover $\gamma$ in $\Omega$ from $\Lambda_\gamma^s$. 

\vspace{5mm}

The following theorems are the main results in this paper. The first one solves the injectivity question relative to the inverse problem for the non-local conductivity equation in any dimension $n\geq 1$. 

\paragraph{\bf Theorem 1.1} \label{1.1} \emph{Let $\Omega\subset\mathbb R^n, \; n\geq 1,$ be a bounded open set, $s\in(0,1)$, and for $j=1,2$ let $\gamma_j: \mathbb R^n \rightarrow \mathbb R$ be such that} 
$$ \left\{\begin{array}{lr}
        \mbox{for some } \underline {\gamma_j}, \overline {\gamma_j} \in \mathbb R, \; \; 0<\underline {\gamma_j} \leq \gamma_j(x) \leq \overline {\gamma_j} < \infty, \; \mbox{for a.e. } x \in \mathbb R^n \\
        \gamma_j^{1/2}(x)-1:= m_j(x) \in W^{2s,n/2s}_c(\Omega)
        \end{array}\right..	 $$

\noindent \emph{Suppose $W_1, W_2\subset \Omega_e$ are open sets, and that the DN maps for the conductivity equations in $\Omega$ relative to $\gamma_1$ and $\gamma_2$ satisfy 
$$\Lambda^s_{{\gamma_1}}[f]|_{W_2}=\Lambda^s_{{\gamma_2}}[f]|_{W_2}, \;\;\;\;\; \forall f\in C^\infty_c(W_1)\;.$$ \noindent Then $\gamma_1 = \gamma_2$. }
\vspace{3mm}

\noindent The second theorem gives a uniqueness result and even a reconstruction procedure for the same inverse problem with a single measurement. 

\paragraph{\bf Theorem 1.2} \label{1.2} \emph{Let $\Omega \subset \mathbb R^n, \; n\geq 1$ be a bounded open set, $s\in (0,1)$, $\epsilon>0$, and let $\gamma: \mathbb R^n\rightarrow \mathbb R$ be such that}  
$$ \left\{\begin{array}{lr}
        \mbox{for some } \underline {\gamma}, \overline {\gamma} \in \mathbb R, \; \; 0<\underline {\gamma} \leq \gamma(x) \leq \overline {\gamma} < \infty, \; \mbox{for a.e. } x \in \mathbb R^n \\
        \gamma^{1/2}(x)-1:= m(x) \in W^{2s+\epsilon,p}_c(\Omega), \;\; \mbox{for } p> n/\epsilon
        \end{array}\right..	 $$

\noindent\emph{Suppose $W_1, W_2 \subset \Omega_e$ are open sets, with $\overline \Omega \cap \overline{W_1} = \emptyset$. Given any fixed function $g\in \tilde H^s(W_1)\setminus \{0\}$, $\gamma$ is uniquely determined and can be reconstructed from the knowledge of $\Lambda^s_\gamma [g]|_{W_2}$.}

\paragraph{\color{black} Remark} {\color{black} In the theorems above we make some regularity assumptions on $m$: namely, it is required to belong to Sobolev spaces of the form $W_c^{k, p}(\Omega)$, which are defined in Section \ref{prel}. Such assumptions are needed in order to be able to apply the previous results (\cite{RS2017}, \cite{GRSU2018}), which are recalled in Section \ref{proof} and constitute the core of the proofs of our theorems.  }
\vspace{3mm}

\noindent A tool that is often used for treating the second order conductivity equation is Liouville's reduction, which consists in rephrasing the problem in terms of the function $w=\gamma^{1/2}u$ and the potential $q=\frac{\Delta\gamma^{1/2}}{\gamma^{1/2}}$. It is easily shown that the resulting equation is $-\Delta w + qw =0$, i.e. Schr\"odinger's equation. The idea behind the proofs of Theorems 1.1 and 1.2 is to use a reduction similar to Liouville's, but suited for a non-local setting: as it will be shown in Section \ref{proof}, the potential will be $q=-\frac{(-\Delta)^sm}{\gamma^{1/2}}$. The problems considered are thus transformed into special cases of inverse problems for the fractional Schr\"odinger equation. These are in turn well understood and dealt with thanks to the previous results (\cite{RS2017}, \cite{GRSU2018}). The key points in these works are the strong uniqueness and approximation results obtained in \cite{DSV2016}. For an overview of the fractional Calder\'on problem and many more references, see the survey \cite{S2017}.

\vspace{5mm}

This paper is organized as follows. Section \ref{intro} is the introduction. Section \ref{prel} is devoted to preliminaries and definitions, including Sobolev spaces and non-local operators. Section \ref{proof} first defines the conductivity equation and the DN map, then proves the main theorems. Section \ref{lim} contains an analysis of the limit case $s\rightarrow 1^-$, which is expected to give the local problem. The last part, Section \ref{model}, is devoted to a simple model for a random walk with long jumps from which the fractional conductivity equation naturally arises.  
\vspace{5mm}

\noindent {\bf Acknowledgement.} This work is part of the PhD research of the author. The author was partially supported by the European Research Council under Horizon 2020 (ERC CoG 770924). The author wishes to dearly thank professor M. Salo for his precious ideas and helpful discussion in the making of this work.   

\section{Preliminaries} \label{prel}

In this section the main function spaces, operators and notations of the paper will be introduced. For the Sobolev spaces, the notation will be the usual one (check, e.g., \cite{GSU2017}). The non-local operators are based on the theoretical framework presented in \cite{DGLZ2012}.

\paragraph{\bf Sobolev spaces} {\color{black} If $k\in \mathbb R$, $p\in(1,\infty)$ and $n\in\mathbb{N}\setminus \{0\}$, the symbols $W^{k,p} = W^{k,p}(\mathbb R^n)$ indicate the usual $L^p$-based Sobolev space. If $\Omega\subset \mathbb R^n$ is an open set, the symbol $W^{k,p}_c(\Omega)$ indicates that subset of $W^{k,p}$ whose elements can be approximated in the Sobolev norm by functions belonging to $C^\infty_c(\Omega)$. }

\noindent In particular, given $s\in (0,1)$ and $n\in\mathbb{N}\setminus \{0\}$, the symbols $H^s = H^s(\mathbb{R}^n) = W^{s,2}(\mathbb{R}^n)$ indicate the standard $L^2$-based Sobolev space with norm $$\|u\|_{H^s(\mathbb{R}^n)} = \|\mathcal{F}^{-1}( \langle\xi\rangle^s \hat u ) \|_{L^2(\mathbb{R}^n)}\;,$$ where $\langle\xi\rangle := (1+|\xi|^2)^{1/2}$. The notation for the Fourier transform is $\hat u(\xi) = \mathcal F u(\xi) = \int_{\mathbb R^n} e^{-ix\cdot\xi} u(x) dx\;.$ If $U, F\subset \mathbb R^n$ are an open and a closed set, define $$H^s(U) = \{ u|_U, u\in H^s(\mathbb R^n) \}\;,$$$$ \tilde H^s(U) = \text{closure of $C^\infty_c(U)$ in $H^s(\mathbb R^n)$}\;, \;\mbox{and}$$ $$ H^s_F(\mathbb R^n) = \{ u\in H^s(\mathbb R^n) : \text{supp}(u) \subset F \}\;. $$

\noindent The set $H^s(U)$ is equipped with the norm $\|u\|_{H^s(U)}= \inf\{ \|w\|_{H^s(\mathbb R^n)} ; w\in H^s(\mathbb R^n), w|_{U}=u \}$. If $U$ is a Lipschitz domain, the Sobolev spaces $\tilde H^s(U)$ and $H^s_{\bar U}(\mathbb R^n)$ can be naturally identified for all real $s$. For more details on this topic, check \cite{GSU2017}.  
\vspace{2mm}

If $U\subset \mathbb R^n$ is a bounded open set and $s\in(0,1)$, let $X= H^s(\mathbb R^n)/\tilde H^s(U)$ be the \emph{abstract trace space}. If $U$ is a Lipschitz domain, $X$ is the quotient $H^s(\mathbb R^n)/H^s_{\bar U}(\mathbb R^n)$, in which two functions $u,v\in H^s(\mathbb R^n)$ are equivalent if and only if $u|_{U_e} = v|_{U_e}$. 

\paragraph{Remark} There exist several other definitions of Sobolev spaces. In fact (\cite{hitch}, prop. 3.4), given $s\in (0,1)$ and an open set $U\subset \mathbb R^n$ whose boundary is regular enough (in the sense of \cite{hitch}, prop. 2.2), $H^s(U)$ might just as well be defined as $$\check H^s(U) = \left\{ u\in L^2(U) : \frac{|u(x)-u(y)|}{|x-y|^{n/2+s}}\in L^2(U^2) \right\}\; ,$$ with the natural norm \begin{equation}\label{eq:sobodef2}\begin{split} \|u\|_{\check H^s(U)} = \left( \|u\|_{L^2(U)}^2 + [u]_{\check H^s(U)}^2 \right)^{1/2} \; ,\\ [u]_{\check H^s(U)} := \left( \int_U\int_U \frac{|u(x)-u(y)|^2}{|x-y|^{n+2s}} dx\,dy \right)^{1/2}\;. \end{split}\end{equation}

\vspace{3mm}
\paragraph{\bf Non-local operators} If $u\in \mathcal S(\mathbb R^n)$, its fractional Laplacian is $$ (-\Delta)^s u(x) := \mathcal C_{n,s} \lim_{\epsilon \rightarrow 0^+}\int_{\mathbb R^n\setminus B_\epsilon (x)} \frac{u(y)-u(x)}{|y-x|^{n+2s}} dy\;,$$ where $\mathcal C_{n,s} := \frac{4^s \Gamma(n/2+s)}{\pi^{n/2}|\Gamma(-s)|}$ is a constant satisfying (see \cite{hitch}) \begin{equation}
\label{eq:Cns, sto1}
\begin{split}
\lim_{s\rightarrow 1^-} \frac{\mathcal C_{n,s}}{s(1-s)} = \frac{4n}{\omega_{n-1}}\;.
\end{split}
\end{equation} This choice assures that the Fourier symbol of the fractional Laplacian is $|\xi|^{2s}$, i.e. the equality $(-\Delta)^s u(x) =\mathcal F^{-1} ( |\xi|^{2s} \hat u (\xi) )$ holds. If $k\in \mathbb R$ and $p\in(1,\infty)$, $(-\Delta)^s$ extends as a bounded map (\cite{Ho90}, Ch. 4 and \cite{Ta96}) \begin{equation}\label{lapmapsto} (-\Delta)^s : W^{k,p}(\mathbb R^n) \rightarrow W^{k-2s,p}(\mathbb R^n)\;.\end{equation}
\noindent {\color{black} For the sake of completeness, it should be added that there exist many equivalent definitions for the fractional Laplacian (\cite{Kw15}). As shown by change of variables in \cite{hitch}, one of them is }
\begin{equation}\label{eq:lapother}\color{black} (-\Delta)^sv(x)=-\frac{\mathcal C_{n,s}}{2} \; PV \int_{\mathbb R^n} \frac{\delta v(x,y)}{|y|^{n+2s}} dy\;, \end{equation}

\noindent {\color{black} which holds if $v$ is a Schwartz function. The symbol $\delta v(x,y)$, which is quite recurrent in this paper, is defined as follows:}
\begin{equation}
\label{eq:delta} \color{black}
\delta v (x,y) := v(x+y)+v(x-y)-2v(x)\;.
\end{equation}

\noindent {\color{black} This way of writing the fractional Laplacian is very useful for removing the singularity at the origin: in fact, if $v$ is a smooth function, by means of a Taylor expansion one gets
$$ \frac{v(x+y) + v(x-y) - 2v(x)}{|y|^{n+2s}} \leq \frac{\|D^2 v\|_{L^\infty}}{|y|^{n+2s-2}}\;, $$

\noindent which is integrable near 0.}
\vspace{3mm}

Motivated by the elementary decomposition $\Delta u = \nabla\cdot(\nabla u)$, the next step will be to define two fractional counterparts of such differential operators, following \cite{DGLZ2012}. These will share the non-local nature of the fractional Laplacian. 
\vspace{3mm}

\noindent Let $u\in C^{\infty}_c(\mathbb R^n)$, $x,y \in \mathbb R^n$. The \emph{fractional gradient of $u$ at points $x$ and $y$} is the vector \begin{equation}\label{eq:graddef}\nabla^s u(x,y) := -\frac{\mathcal C_{n,s}^{1/2}}{\sqrt 2} \frac{u(y)-u(x)}{|y-x|^{n/2 +s+1}}(y-x)\;.\end{equation}

\noindent Using the result (prop. 3.6, \cite{hitch}), formula (\ref{eq:sobodef2}) and the fact that $0\leq |\xi|/\langle\xi\rangle \leq 1$, it is easy to find the following inequality:  
\begin{equation}
\label{eq:formula61/2}
\begin{split}
\|\nabla^s u\|^2_{L^2(\mathbb R^{2n})} & = \frac{C_{n,s}}{2} \int_{\mathbb R^n}\int_{\mathbb R^n} \frac{|u(x)-u(y)|^2}{|x-y|^{n+2s}} dx\,dy = \frac{C_{n,s}}{2} [u]^2_{\check H^s(\mathbb R^n)} \\ & = \| (-\Delta)^{s/2} u \|^2_{L^2(\mathbb R^n)} = \| |\xi|^s \hat u(\xi) \|^2_{L^2(\mathbb R^n)} = \left\| \frac{|\xi|^s}{\langle\xi\rangle^s} \langle\xi\rangle^s\hat u(\xi) \right\|^2_{L^2(\mathbb R^n)} \\ & \leq \| \langle\xi\rangle^s\hat u(\xi) \|^2_{L^2(\mathbb R^n)} = \|u\|^2_{H^s(\mathbb R^n)}\;.
\end{split}
\end{equation}

\noindent Thus the linear operator $\nabla^s$ maps $C^\infty_c(\mathbb R^n)$ into $L^2(\mathbb R^{2n})$. What is more, since $C^\infty_c(\mathbb R^n)$ is dense in $H^s(\mathbb R^n)$ one can extend $\nabla^s$ so that it is defined in $H^s(\mathbb R^n)$ and formula (\ref{eq:formula61/2}) still holds. 

\noindent The next operator is defined by duality. Let $u\in H^s(\mathbb R^n)$, $v\in L^2(\mathbb R^{2n})$; the \emph{fractional divergence} is that operator $(\nabla\cdot)^s : L^2(\mathbb R^{2n}) \rightarrow H^{-s}(\mathbb R^n)$ such that the following formula holds:
\begin{equation}
\label{eq:defdiv}
\langle (\nabla\cdot)^s v,u \rangle_{L^2(\mathbb R^{n})} = \langle v,\nabla^s u \rangle_{L^2(\mathbb R^{2n})}  \;.
\end{equation}

The next simple lemma allows the composition of the fractional divergence and its adjoint into the fractional Laplacian.

\paragraph{\bf Lemma 2.1} \label{2.0} \emph{\color{black} Let $u\in H^s(\mathbb R^n)$. Then the equality $(\nabla\cdot)^s(\nabla^su)(x) = (-\Delta)^su(x)$ holds in weak sense, with $(\nabla\cdot)^s(\nabla^su) \in H^{-s}(\mathbb R^n)$.} 

\paragraph{Proof} Let $u$, $\phi \in H^s(\mathbb R^n)$, and by density for all $i\in \mathbb N$ let $u_i, \phi_i$ be smooth, compactly supported functions such that $\|u-u_i\|_{H^s(\mathbb R^n)} \leq 1/i$ and $\|\phi-\phi_i\|_{H^s(\mathbb R^n)} \leq 1/i$. By Cauchy-Schwartz inequality and formula (\ref{eq:formula61/2}), it is seen that
\begin{equation*}
\langle \nabla^su, \nabla^s\phi \rangle = \lim_{i\rightarrow\infty} \left( \langle \nabla^s(u-u_i), \nabla^s\phi \rangle + \langle \nabla^su_i, \nabla^s\phi \rangle \right) = \lim_{i\rightarrow\infty} \langle \nabla^su_i, \nabla^s\phi \rangle\;,
\end{equation*}
\noindent and thus $\langle \nabla^su, \nabla^s\phi \rangle = \lim_{i\rightarrow\infty} \langle \nabla^su_i, \nabla^s\phi_i \rangle$. Now compute
\begin{equation*} \color{black}
\begin{split} 
\langle \nabla^su_i,& \nabla^s\phi_i \rangle  = \frac{C_{n,s}}{2} \int_{\mathbb R^n}\int_{\mathbb R^n} \frac{u_i(y)-u_i(x)}{|y-x|^{n+2s}} (\phi_i(y)-\phi_i(x))\; dydx \\ & = \frac{C_{n,s}}{2} \int_{\mathbb R^n}\int_{\mathbb R^n} \frac{u_i(x\pm z)-u_i(x)}{|z|^{n+2s}} (\phi_i(x\pm z)-\phi_i(x))\; dzdx \\ & = \frac{C_{n,s}}{4} \int_{\mathbb R^n}\int_{\mathbb R^n} \frac{1}{|z|^{n+2s}}\Big\{ -\phi_i(x) \delta u_i(x,z) +(u_i\phi_i)(x+z) + (u_i\phi_i)(x-z)  \Big.\\ & \Big. \;\;\;\;\;\;\;\;\;\;\;\;\;\;\;\;\;\;\;\;\;\;\;\;\;\;\;\;\;\;\;\;\;\;\;\;\;\;\,\; -u_i(x) (\phi_i(x+z) + \phi_i(x-z) )\Big\} \;dzdx\;.
\end{split}
\end{equation*}
\noindent By adding and subtracting the term $2 (u_i\phi_i)(x)$ we then get
\begin{equation*}
\begin{split}
\langle \nabla^su_i,& \nabla^s\phi_i \rangle  = \frac{C_{n,s}}{4} \int_{\mathbb R^n}\int_{\mathbb R^n} \frac{  -\phi_i(x) \delta u_i(x,z)  + \delta (u_i\phi_i)(x,z) - u_i(x) \delta \phi_i (x,z)  }{|z|^{n+2s}}\; dzdx. 
\end{split}
\end{equation*}

\noindent This integral can be split in three parts, which are all well defined because of the above consideration about the removal of the singularity at the origin:
\begin{equation*}
\begin{split}
\langle \nabla^su_i, \nabla^s\phi_i \rangle  & = \frac{1}{2} \Big( \langle \phi_i, (-\Delta)^s u_i \rangle - \langle 1, (-\Delta)^s(u_i\phi_i)\rangle + \langle u_i, (-\Delta)^s\phi_i\rangle\Big) \\ & = \langle \phi_i, (-\Delta)^su_i \rangle\;.
\end{split}
\end{equation*}

\noindent The last line follows from the fact that $u_i, \phi_i \in C^\infty_c(\mathbb R^n)$, which means that the first and last terms are equal. {\color{black} Moreover, the second term vanishes because, by Fubini's theorem, \begin{equation*}\begin{split} \langle 1, (-\Delta)^s(u_i\phi_i)\rangle & = -\frac{ \mathcal C_{n,s}}{2} \int_{\mathbb R^n} \int_{\mathbb R^n} \frac{ (u_i\phi_i)(x+y) +(u_i\phi_i)(x-y) - 2(u_i\phi_i)(x)}{|y|^{n+2s}} \, dydx \\ & = -\frac{ \mathcal C_{n,s}}{2} \int_{\mathbb R^n} \frac{1}{|y|^{n+2s}}\int_{\mathbb R^n} ( (u_i\phi_i)(x+y) +(u_i\phi_i)(x-y) - 2(u_i\phi_i)(x)) \, dxdy\;,\end{split}\end{equation*} \noindent and the integral in $dx$ is of course independent of $y$ and equal to $0$.} Therefore $\langle \nabla^su_i, \nabla^s\phi_i \rangle = \langle (-\Delta)^su_i, \phi_i \rangle$, and eventually
\begin{equation*}
\begin{split}
\langle (\nabla\cdot)^s(\nabla^s)u,&\phi \rangle  := \langle \nabla^su, \nabla^s\phi \rangle = \lim_{i\rightarrow\infty} \langle \nabla^su_i, \nabla^s\phi_i \rangle = \lim_{i\rightarrow\infty} \langle (-\Delta)^su_i, \phi_i \rangle \\ & = \lim_{i\rightarrow\infty} \Big( \langle (-\Delta)^s(u_i-u), \phi_i \rangle + \langle (-\Delta)^su, \phi_i-\phi \rangle \Big) + \langle (-\Delta)^su, \phi \rangle \\ & =  \langle (-\Delta)^su, \phi \rangle \;,
\end{split}
\end{equation*}
\noindent just as wanted. Notice that the limit vanishes because $\|(-\Delta)^sw\|_{H^{-s}} \leq \|w\|_{H^s}$. This proves the first statement; the second one now follows from the previous remark about the extensions of the fractional Laplacian. \qed

{\color{black} \paragraph{Remark} \noindent $\nabla^s$ and $(\nabla\cdot)^s$ can be respectively identified with the operators $\mathcal D^*$ and $\mathcal D$ from \cite{DGLZ2012}, where the antisymmetric vector mapping $\alpha(x,y): \mathbb R^{2n}\rightarrow \mathbb R^n$ is chosen as \begin{equation}\label{eq:alphadef}\alpha(x,y)= \frac{\mathcal C_{n,s}^{1/2}}{\sqrt 2} \frac{y-x}{|y-x|^{n/2 +s+1}}\;.\end{equation} \noindent The choice of $\alpha$ comes from the fact that we want to have $(\nabla\cdot)^s(\nabla^su)=(-\Delta)^su$, which at least for $u\in \mathcal S$ means $$ (-\Delta)^su(x) = 2\int_{\mathbb R^n} (u(x)-u(y))|\alpha(x,y)|^2 dy\;. $$ \noindent Thus the most natural choice would be to have $|\alpha(x,y)| = \frac{\mathcal C_{n,s}^{1/2}}{\sqrt{2} |y-x|^{n/2+s}}$, which motivates our choice of $\alpha$. In this case we also have, for $u\in C^\infty_c(\mathbb R^n)$, \begin{equation}\label{modgrad} |\nabla^su|= \frac{\mathcal C_{n,s}^{1/2}}{\sqrt{2}}\frac{|u(y)-u(x)|}{|y-x|^{n/2+s}}\;. \end{equation} \noindent Anyway, different choices of $\alpha$ could in principle be considered.}

\section{Main results} \label{proof}

\paragraph{\bf Non-local conductivity equation} Let $\Omega \subset \mathbb R^n$ be an open set; we call $\Omega_e = \mathbb R^n\setminus \overline\Omega$ the \emph{exterior domain}.

{\color{black} Let $\gamma : \mathbb R^n \rightarrow \mathbb R$ be a measurable function such that there exist $\underline \gamma, \overline \gamma \in \mathbb R$ such that $ 0<\underline \gamma \leq \gamma(x) \leq \overline \gamma < \infty$ for all $x \in \mathbb R^n$,} and let $m(x):= \gamma^{1/2}(x)-1$ belong to $W^{2s,n/2s}_c(\Omega)$. The assumptions for the \emph{conductivity} $\gamma$ are similar to the ones that are typically made in the second order case. The values of $\gamma(x)$ for $x\in$ supp $m$ represent the conductivity in the object of study. Outside of this region $\gamma(x)\equiv 1$, because the electrical properties of the surroundings are thought of as constant. 
\vspace{2mm}

Let $\Theta : \mathbb R^{2n} \rightarrow \mathbb R^{2n}$ be the variable matrix $\Theta(x,y) := \gamma(x)^{1/2}\gamma(y)^{1/2} \text{Id}$. The \emph{interaction matrix} $\Theta$ represents how readily an electron will jump from $x$ to $y$. We assume the material to be isotropic, meaning that the interaction does not depend on direction; therefore, $\Theta(x,y)$ is a symmetrical scalar multiple of the identity matrix.  

{\color{black}\paragraph{Remark} According to formula \eqref{lapmapsto}, it makes sense to compute $(-\Delta)^sm$, and it belongs to $W^{0,n/2s}(\mathbb R^n) = L^{n/2s}(\mathbb R^n)$.}

\vspace{5mm}
By using the boundedness of $\gamma$ and Lemma 2.1 it is seen that if $u\in H^s(\mathbb R^n)$, then $\Theta\cdot\nabla^su\in L^2(\mathbb R^{2n})$:
\begin{equation*}
\|\Theta\cdot\nabla^su\|^2_{L^2(\mathbb R^{2n})} = \int_{\mathbb R^{2n}} \gamma(x)\gamma(y) \nabla^su\cdot\nabla^su \;dx\;dy \leq \overline{\gamma}^2 \|\nabla^su\|^2_{L^2(\mathbb R^{2n})} < \infty. 
\end{equation*}

\noindent Let $u\in H^s(\mathbb R^n)$. The \emph{non-local conductivity operator} is $\mathrm C_\gamma^s u := (\nabla\cdot)^s(\Theta\cdot\nabla^s u)$, while the \emph{non-local conductivity equation} is the statement $\mathrm C_\gamma^s u=0$ in $\Omega$.

The next theorem reduces the conductivity equation to Schr\"odinger's.

\paragraph{\bf Theorem 3.1} \label{3.2} \emph{Let $u\in H^s(\mathbb R^n)$, $g\in H^s(\Omega_e)$, $w=\gamma^{1/2}u$, $f=\gamma^{1/2}g$ and $q=-\frac{(-\Delta)^sm}{\gamma^{1/2}}$. $u$ solves the conductivity equation with exterior value $g$ if and only if $w$ solves Schr\"odinger's equation with exterior value $f$, that is}
\begin{equation*} 
\left\{\begin{array}{lr}
        (\nabla\cdot)^s(\Theta\cdot\nabla^s u)=0 & \text{in } \Omega\,\;\\
        u=g & \text{in } \Omega_e
        \end{array}\right. \;\;\Leftrightarrow\;\; 
		\left\{\begin{array}{lr}
        \Big( (-\Delta)^s +q \Big)w=0 & \text{in } \Omega\,\;\\
        w=f & \text{in } \Omega_e
        \end{array}\right. \;.		 
\end{equation*}

\noindent \emph{Moreover, the following formula holds for all $w\in H^s(\mathbb R^n)$:}
$$ \mathrm C^s_\gamma (\gamma^{-1/2}w) = \gamma^{1/2} ((-\Delta)^s+q)w\,. $$

\paragraph{Proof} Start by observing that $m$ is a Fourier multiplier on $H^s$, because we have the embedding $(W^{2s,n/2s}\cap L^\infty) \times H^s \hookrightarrow H^s$ (check Lemma 6, \cite{BM01}). This of course means that also $\gamma^{1/2}=1+m$ is a Fourier multiplier on $H^s$, which in turn implies that $w\in H^s$ and $f\in H^s(\Omega_e)$. Moreover, the computation $$ qw = -\frac{(-\Delta)^sm}{\gamma^{1/2}} \,\gamma^{1/2}u = -u(-\Delta)^sm $$ and the observation that, by Theorem 6.1 in \cite{BH2017} and Sobolev embedding theorem, $$L^{n/2s} \times H^s \hookrightarrow L^{2n/(n+2s)} \hookrightarrow H^{-s}$$ imply that $((-\Delta)^s + q)w \in H^{-s}$.
\vspace{2mm}

\noindent Our proof will be very similar to the one of the previous Lemma 2.1. Take $\phi \in H^s$, and for all $i\in\mathbb N$ let $\phi_i, u_i\in C^\infty_c(\mathbb R^n)$ be such that $\|\phi-\phi_i\|_{H^s}<1/i$ and $\|u-u_i\|_{H^s}<1/i$. By definition, Cauchy-Schwartz inequality and formula (\ref{eq:formula61/2}) we get
\begin{equation}
\label{eq:new3.1.1}
\begin{split} 
\langle \mathrm C_\gamma^s u, \phi \rangle & = \langle (\nabla\cdot)^s(\Theta\cdot\nabla^su), \phi \rangle = \langle \Theta\cdot\nabla^su,\nabla^s\phi \rangle \\ & = \lim_{i\rightarrow \infty} \Big( \langle \Theta\cdot\nabla^su,\nabla^s\phi_i \rangle  + \langle \Theta\cdot\nabla^su,\nabla^s(\phi-\phi_i) \rangle\Big) \\ & = \lim_{i\rightarrow \infty} \langle \Theta\cdot\nabla^su,\nabla^s\phi_i \rangle = \lim_{i\rightarrow \infty} \langle \Theta\cdot\nabla^su_i,\nabla^s\phi_i \rangle\;.
\end{split}
\end{equation}

\noindent By change of variables,
\begin{equation*}
\begin{split}
\langle \Theta&\cdot\nabla^su_i,\nabla^s\phi_i \rangle  = \frac{\mathcal C_{n,s}}{2} \int_{\mathbb R^n}\int_{\mathbb R^n} \gamma(x)^{1/2}\gamma(y)^{1/2}\, \frac{(u_i(y)-u_i(x))\,(\phi_i(y)-\phi_i(x))}{|y-x|^{n +2s}} dydx \\ & = \frac{\mathcal C_{n,s}}{2} \int_{\mathbb R^n}\int_{\mathbb R^n} \gamma(x)^{1/2}\gamma(x\pm z)^{1/2}\, \frac{(u_i(x\pm z)-u_i(x))\,(\phi_i(x\pm z)-\phi_i(x))}{|z|^{n +2s}} dzdx \\ & = \frac{\mathcal C_{n,s}}{4} \int_{\mathbb R^n}\int_{\mathbb R^n} \Big\{ \gamma(x)^{1/2}\gamma(x+ z)^{1/2}\, \frac{(u_i(x+ z)-u_i(x))\,(\phi_i(x+ z)-\phi_i(x))}{|z|^{n +2s}} \Big. \\ & \;\;\;\;\;\;\;\;\;\;\;\;\;\;\;\;\;\;\;\;\;\;\;\;\;\; \Big.   + \gamma(x)^{1/2}\gamma(x- z)^{1/2}\, \frac{(u_i(x- z)-u_i(x))\,(\phi_i(x- z)-\phi_i(x))}{|z|^{n +2s}}\Big\}\,dzdx \;.
\end{split}
\end{equation*}
\noindent Now consider the integrand function. By defining $w_i := \gamma^{1/2}u_i$ it can be rewritten as
\begin{equation*}
\begin{split}
 \frac{\gamma(x)^{1/2}}{|z|^{n+2s}}\Big\{&  -\phi_i(x)\left( w_i(x+z) + w_i(x-z) - u_i(x)( \gamma^{1/2}(x+z) + \gamma^{1/2}(x-z) ) \right)  +   \Big.\\ & \Big.     
       (w_i\phi_i)(x+z) + (w_i\phi_i)(x-z) - u_i(x)\left( (\gamma^{1/2}\phi_i)(x+z) + (\gamma^{1/2}\phi_i)(x-z) \right)                      \Big\}\,, 
\end{split}
\end{equation*}

\noindent so that, if we add and subtract the term $2w_i(x)$ from the first line and the term $ 2(w_i\phi_i)(x)$ from the second one, by formula (\ref{eq:delta}) we get
\begin{equation*}
\begin{split}
 \frac{\gamma(x)^{1/2}}{|z|^{n+2s}}\Big\{ \delta(w_i\phi_i)(x,z) - u_i(x) \delta (\gamma^{1/2}\phi_i)(x,z)  -\phi_i(x)\left( \delta w_i(x,z) - u_i(x) \delta (\gamma^{1/2}-1)(x,z)  \right)\Big\}\,.
\end{split}
\end{equation*}

\noindent Therefore
\begin{equation*}
\begin{split}
\langle \Theta\cdot\nabla^su_i,\nabla^s\phi_i \rangle  = \frac{\mathcal C_{n,s}}{4} \int_{\mathbb R^n}\int_{\mathbb R^n} \frac{\gamma(x)^{1/2}}{|z|^{n+2s}}\Big\{& \delta(w_i\phi_i)(x,z) - u_i(x) \delta (\gamma^{1/2}\phi_i)(x,z)  \Big.\\ & \Big.-\phi_i(x)\left( \delta w_i(x,z) - u_i(x) \delta (\gamma^{1/2}-1)(x,z)  \right)\Big\}\,,
\end{split}
\end{equation*}

\noindent and the interior integral can be split in the following four parts by Lemma 2.1, since the $\delta$'s make each of them integrable at the origin:
\begin{equation*}
\begin{split}
\langle \Theta\cdot\nabla^su_i,\nabla^s\phi_i \rangle  & = \frac{1}{2} \int_{\mathbb R^n}\Big\{ -\gamma^{1/2} (-\Delta)^s(w_i\phi_i)  +  w_i(-\Delta)^s(\gamma^{1/2}\phi_i)  \Big.\\ &\;\;\;\;\;\;\;\;\;\;\;\;\;\;\,\;\;\Big. + \phi_i\gamma^{1/2} (-\Delta)^sw_i -  \phi_i\gamma^{1/2} u_i(-\Delta)^s(\gamma^{1/2}-1)   \Big\} \\ & = \frac{1}{2} \int_{\mathbb R^n}\Big\{ (1-\gamma^{1/2}) (-\Delta)^s(w_i\phi_i)  +  w_i(-\Delta)^s(\gamma^{1/2}\phi_i)  \Big.\\ & \;\;\;\;\;\;\;\;\;\;\;\;\;\;\,\;\;\Big. + \phi_i\gamma^{1/2} (-\Delta)^sw_i -  \phi_i\gamma^{1/2} w_i \frac{(-\Delta)^s(\gamma^{1/2}-1)}{\gamma^{1/2}}   \Big\}\,.
\end{split}
\end{equation*}

\noindent In the last line, we have added the term $\frac{1}{2}\int_{\mathbb R^n} (-\Delta)^s(w_i\phi_i)$, which equals $0$. Now by the first part of the proof we can compute
\begin{equation*}\color{black}\begin{split} \langle \Theta\cdot\nabla^su_i,\nabla^s\phi_i \rangle & = \frac{\langle \gamma^{1/2}\phi_i , ((-\Delta)^s+q)w_i \rangle}{2} + \frac{\langle -(\gamma^{1/2}-1), (-\Delta)^s (w_i\phi_i) \rangle + \langle w_i, (-\Delta)^s (\gamma^{1/2}\phi_i) \rangle}{2} \\ & = \frac{\langle \gamma^{1/2}\phi_i , ((-\Delta)^s+q)w_i \rangle}{2} + \frac{\langle -((-\Delta)^sm) u_i, \gamma^{1/2}\phi_i \rangle + \langle (-\Delta)^sw_i, \gamma^{1/2}\phi_i \rangle}{2} \\ & = \frac{\langle \gamma^{1/2}\phi_i , ((-\Delta)^s+q)w_i +qw_i + (-\Delta)^sw_i\rangle}{2} \\ & = \langle \gamma^{1/2}\phi_i , ((-\Delta)^s+q)w_i \rangle \;.\end{split}\end{equation*}
\noindent Eventually, by using this and (\ref{eq:new3.1.1}),
$$ \langle \mathrm C_\gamma^s u, \phi \rangle = \lim_{i\rightarrow\infty} \langle \gamma^{1/2}\phi_i , ((-\Delta)^s+q)w_i \rangle = \langle \phi , \gamma^{1/2}((-\Delta)^s+q)w \rangle\;.$$
\noindent This last step holds true because
$$ \lim_{i\rightarrow\infty} |\langle \gamma^{1/2}(\phi_i-\phi), ((-\Delta)^s+q)w_i \rangle| \leq c \lim_{i\rightarrow\infty} \|\phi_i - \phi\|_{H^s} \| ((-\Delta)^s+q)w_i \|_{H^{-s}} = 0$$ 
\noindent and
\begin{equation*}\begin{split} \lim_{i\rightarrow\infty} |\langle \gamma^{1/2}\phi, ((-\Delta)^s+q)(w_i-w) &\rangle|  \leq c \|\phi\|_{H^s} \lim_{i\rightarrow\infty} \| ((-\Delta)^s+q)(w_i-w) \|_{H^{-s}} \\ & \leq c \|\phi\|_{H^s} (1+\|q\|_{L^{n/2s}}) \lim_{i\rightarrow\infty}\|w_i-w\|_{H^s} =0\;.\end{split}\end{equation*}\qed

\paragraph{\bf Bilinear form} Let $s\in(0,1),\; u,v\in H^s(\mathbb R^n)$, and define the \emph{bilinear form} $B^s_\gamma : H^s \times H^s \rightarrow \mathbb R$ as follows
\begin{equation}
\label{eq:bilindef}
B^s_\gamma [u,v] = \int_{\mathbb R^n}\int_{\mathbb R^n} \nabla^s v \cdot (\Theta \cdot \nabla^s u) dy \, dx\;.
\end{equation}

$B^s_\gamma$ is a useful instrument to show the well-posedness of the direct problem for the conductivity equation. In \cite{DGLZ2012}, Theorem 4.9, it is proved that for all $F \in (\tilde H^s(\Omega))^*$ there exists a unique solution $u_F \in \tilde H^s(\Omega)$ to $B^s_\gamma[u,v]=F(v), \; \forall v\in \tilde H^s(\Omega)$. This is equivalent to saying that for all $F \in (\tilde H^s(\Omega))^*$ there exists one and only one $u_F\in H^s(\Omega)$ such that $\mathrm C^s_\gamma u =F$ in $\Omega$, $u_F|_{\Omega_e}=0$. To treat the case of non-zero exterior value, suppose $f\in H^s(\mathbb R^n)$ and let $u_f = \overline u +f$, where $\overline u \in H^s(\Omega)$ is the unique solution to the problem  
$$ \left\{\begin{array}{lr}
        \mathrm C^s_\gamma u=F - B^s_\gamma[f,\cdot] & \text{in } \Omega\,\;\\
        u=0 & \text{in } \Omega_e
        \end{array}\right.. \;\;\;\mbox{ Then }\;\;\; \left\{\begin{array}{lr}
        \mathrm C^s_\gamma u=F & \text{in } \Omega\,\;\\
        u=f & \text{in } \Omega_e
        \end{array}\right.$$

\noindent has $u_f \in H^s(\mathbb R^n)$ as its unique solution. Moreover, it follows from \cite{GSU2017} that \begin{equation}\label{eq:estidirect}\|u_f\|_{H^s(\mathbb R^n)} \leq c(\|F\|_{(\tilde H^s(\Omega))^*} + \|f\|_{H^s(\mathbb R^n)})\;.\end{equation} 

The next lemma collects some properties of $B^s_\gamma$.

\paragraph{\bf Lemma 3.2} \label{3.4} \emph{Let $v,w\in H^s(\mathbb R^n)$, $f,g\in H^s(\Omega_e)$ and $u_f, u_g\in H^s(\mathbb R^n)$ be such that $\mathrm C^s_\gamma u_f =\mathrm C^s_\gamma u_g = 0$ in $\Omega$, $u_f|_{\Omega_e} = f$ and $u_g|_{\Omega_e} = g$. Then}
\begin{enumerate}
\item $B^s_\gamma[v,w] = B^s_\gamma[w,v]\;$ (symmetry),
\item $|B^s_\gamma[v,w]| \leq \overline\gamma\|v\|_{H^s(\mathbb R^n)}\|w\|_{H^s(\mathbb R^n)}\;$,
\item $B^s_\gamma[u_f,e_g] = B^s_\gamma[u_g,e_f]\;$,
\end{enumerate}
\noindent \emph{where $e_g, e_f \in H^s(\mathbb R^n)$ are extensions of $g,f$ respectively.}

\paragraph{Proof} Symmetry is showed by using (\ref{eq:graddef}) in (\ref{eq:bilindef}), $$B^s_\gamma [v,w] = \frac{\mathcal C_{n,s}}{2} \int_{\mathbb R^n}\int_{\mathbb R^n} \gamma(x)^{1/2}\gamma(y)^{1/2}\frac{(v(y)-v(x))\,(w(y)-w(x))}{|y-x|^{n+2s}}\, dy\, dx\;.$$
\noindent For the second point, using H\"older's inequality and the known estimate for the $L^2$ norm of the fractional gradient, \begin{equation*}\begin{split} |B^s_\gamma[v,w]| & \leq \|\nabla^s v\|_{L^2(\mathbb R^{2n})}\|\Theta \cdot \nabla^s w\|_{L^2(\mathbb R^{2n})} \leq \overline\gamma \|\nabla^s v\|_{L^2(\mathbb R^{2n})}\|\nabla^s w\|_{L^2(\mathbb R^{2n})} \\ & \leq \overline\gamma \|v\|_{H^s(\mathbb R^n)}\|w\|_{H^s(\mathbb R^n)}\;. \end{split}\end{equation*} \noindent In order to prove the last point, use the definition of fractional divergence (\ref{eq:defdiv}) $$ B^s_\gamma [u_f,u_g] = \int_{\mathbb R^n}\int_{\mathbb R^n} \nabla^s u_g \cdot (\Theta \cdot \nabla^s u_f) dy \, dx = \int_{\mathbb R^n} u_g \mathrm C_\gamma^s u_f \,dx \;,$$ \noindent then observe that $C_\gamma^s u_f = 0$ in $\Omega$ and $u_g=g$ in $\Omega_e$, so that actually $$ B^s_\gamma [u_f,u_g] = \int_{\Omega_e} u_g \mathrm C_\gamma^s u_f \,dx = \int_{\Omega_e} g \mathrm C_\gamma^s u_f \,dx = \int_{\mathbb R^n} e_g \mathrm C_\gamma^s u_f \,dx = B_\gamma^s [u_f,e_g]\; .$$ \noindent This completes the proof, since by symmetry $$ B_\gamma^s [u_f,e_g] = B^s_\gamma [u_f,u_g] = B^s_\gamma [u_g,u_f] = B^s_\gamma [u_g,e_f] \;.$$\qed

\vspace{3mm}
\paragraph{\bf DN map} The main use of the bilinear form in this paper is the definition of the DN map. In the case of the fractional Calder\'on problem for the Schr\"odinger equation with an unknown potential $q$, such map is $\Lambda_q: X\rightarrow X^*$, $$\Lambda_q[f]([v]) = \int_{\mathbb R^n} v(-\Delta)^sw_f \, dx + \int_\Omega qvw_f\, dx\;,$$ as defined in \cite{GSU2017}. In the above formula, $f,v\in H^s(\mathbb R^n)$ and $w_f\in H^s(\mathbb R^n)$ is the unique solution to $(-\Delta)^sw +qw=0$ in $\Omega$ with $w-f\in \tilde H^s(\Omega)$.  

\paragraph{\bf Lemma 3.3} \emph{There exists a bounded, linear, self-adjoint map $\Lambda_\gamma^s : X\rightarrow X^*$ defined by $$\langle \Lambda_\gamma^s [f],[g] \rangle = B^s_\gamma[u_f,g], \;\;\;\;\;\;\; \forall f,g\in H^s(\mathbb R^n)\; ,$$ \noindent where $X$ is the abstract quotient space $H^s(\mathbb R^n)/\tilde H^s(\Omega)$ and $u_f\in H^s(\mathbb R^n)$ solves $\mathrm C^s_\gamma u = 0$ in $\Omega$ with $u-f\in \tilde H^s(\Omega)$.}
\vspace{3mm}

\paragraph{Proof} The DN map needs to be well defined, that is for all $\phi,\psi\in \tilde H^s(\Omega)$ and $f,g\in H^s(\mathbb R^n)$ the equality $B^s_\gamma[u_f,g] = B^s_\gamma[u_{f+\phi},g+\psi]$ must hold. By Lemma 3.2,
\begin{equation*}
\begin{split}
B^s_\gamma[u_{f+\phi},g+\psi] & = B^s_\gamma[u_{f+\phi},g] + B^s_\gamma[u_{f+\phi},\psi]  = B^s_\gamma[f+\phi,u_g] + \int \psi \mathrm C_\gamma^s u_{f+\phi} \,dx \\ & = B^s_\gamma[f,u_g] + B^s_\gamma[\phi,u_g] = B^s_\gamma[u_f,g] + \int \phi \mathrm C_\gamma^s u_g \,dx = B^s_\gamma[u_f,g]\;,
\end{split}
\end{equation*}
\noindent since $u_{f+\phi}, u_g$ are solutions to the conductivity equation, and $\phi,\psi$ are supported in $\Omega$. The boundedness of $\Lambda_\gamma^s$ follows from the second point of Lemma 3.2 and equation (\ref{eq:estidirect}). In fact, 
\begin{equation*}\begin{split} |\langle \Lambda_\gamma^s [f],[g] \rangle| & = |B^s_\gamma[u_f,g]| \leq c\|u_f\|_{H^s(\mathbb R^n)}\|g\|_{H^s(\mathbb R^n)} \\ & \leq c\|f\|_{H^s(\mathbb R^n)}\|g\|_{H^s(\mathbb R^n)}\;,\;\;\;\;\;\;\;\;\;\;\forall f\in [f],\;\forall g\in [g]\;,\end{split}\end{equation*}  \noindent which implies $$|\langle \Lambda_\gamma^s [f],[g] \rangle| \leq c \inf_{f\in [f]}\|f\|_{H^s(\mathbb R^n)}\inf_{g\in [g]}\|g\|_{H^s(\mathbb R^n)}= c\|[f]\|_X\|[g]\|_X\;.$$ 

\noindent Self-adjointness is trivial, in light of point (3) of Lemma 3.2 : $$ \langle \Lambda_\gamma^s [f],[g] \rangle = B^s_\gamma [u_f,g] = B^s_\gamma[u_g,f] = \langle \Lambda_\gamma^s [g],[f] \rangle = \langle [f],\Lambda_\gamma^s [g] \rangle\;. $$
\hfill \qed
\vspace{3mm} 

\paragraph{\bf Lemma 3.4} \label{3.6} \emph{Let $f,v\in H^s(\mathbb R^n)$ be such that supp$(f)$, supp$(v) \subset \Omega_e$. The DN maps for the conductivity equation $\Lambda^s_\gamma$ and for the corresponding Schr\"odinger equation $\Lambda_{q_\gamma}$ satisfy} $$\Lambda_{q_\gamma}[f]([v]) -\Lambda^s_\gamma [f]([v]) = \int_{\Omega_e}fv(-\Delta)^sm \, dx\;.$$

\paragraph{Proof} First of all observe that we have $\gamma^{1/2}f = f$ and $\gamma^{1/2}v=v$, because supp$(f) \cap$supp$(m)=\emptyset$ and supp$(v) \cap$supp$(m)=\emptyset$. With this in mind and making use of Theorem 3.1 it is easy to compute
\begin{equation*}
\begin{split}
\Lambda^s_\gamma[f]([v]) & = B^s_\gamma [u_f, v] = \int_{\mathbb R^n}\int_{\mathbb R^n} \nabla^sv \cdot (\Theta\cdot\nabla^s u_f)\,dydx \\ & = \int_{\mathbb R^n} v \mathrm C^s_\gamma u_f \, dx = \int_{\mathbb R^n} v\gamma^{1/2}\, ((-\Delta)^s+q_\gamma) w_f \, dx \\ & = \int_{\mathbb R^n} \gamma^{1/2}v(-\Delta)^sw_f\, dx + \int_{\mathbb R^n} \gamma^{1/2}vq_\gamma w_f\, dx \\ & = \int_{\mathbb R^n} v(-\Delta)^sw_f\, dx - \int_{\Omega_e} vf(-\Delta)^sm\, dx\;.
\end{split}
\end{equation*} 
\noindent Moreover, recalling the assumptions about the supports,
$$ \Lambda_{q_\gamma}[f]([v]) = \int_{\mathbb R^n} v(-\Delta)^sw_f\, dx + \int_{\Omega} q_\gamma vw_f\, dx = \int_{\mathbb R^n} v(-\Delta)^sw_f\, dx \;. $$
\noindent The statement of the Lemma is thus proved by taking the difference of the last two formulas.  \hfill \qed

\vspace{2mm}
The definition of the DN map given above, which is abstract in nature, lets us formulate and solve the inverse problems completely. Nonetheless, in the next theorem we will give a more concrete definition of the DN map under stronger assumptions. 

\paragraph{\bf Theorem 3.5} \label{3.5} \emph{Let $\Omega$ be a bounded open set with $C^\infty$ boundary, let $s\in(0,1)$ and let $\gamma^{1/2}=1+m$, with $m\in C^\infty_c(\Omega)$ and $0<\underline \gamma \leq \gamma(x)$, for all $x\in \mathbb R^n$. For any $\beta\geq 0$ such that $\beta\in(s-1/2, 1/2)$ the restriction of $\Lambda^s_\gamma$ to $H^{s+\beta}(\Omega_e)$ is the map} $$ \Lambda^s_\gamma : H^{s+\beta}(\Omega_e) \rightarrow H^{-s+\beta}(\Omega_e), \;\; \Lambda^s_\gamma f = \mathrm C^s_\gamma u_f|_{\Omega_e}\;, $$
\noindent \emph{where $u_f \in H^{s+\beta}(\mathbb R^n)$ solves $\mathrm C^s_\gamma u=0$ in $\Omega$ with $u|_{\Omega_e} =f$, $f\in H^{s+\beta}(\Omega_e)$.}

\paragraph{Proof} Start by observing that the embedding $H^a \times H^c \hookrightarrow H^c$ can be made to work for any $c\in\mathbb R$, if $a$ is taken accordingly large enough: in the case $c<0$, use Theorem 8.1 from \cite{BH2017} with $a>n/2$, while if $c\geq 0$ use Theorem 7.3 with $a>\max\{ n/2,c  \}$. Since now $m\in C^\infty_c(\Omega) \subset H^a(\mathbb R^n)$ for all $a\geq 0$, and consequently $(-\Delta)^s m \in H^{a-2s}$ for all $a\geq 0$, we have that $h\in H^c$ implies $mh, (-\Delta)^sm\, h \in H^c$. It also easily follows that $\gamma^{1/2}h, \gamma^{-1/2} h \in H^c$. 
\vspace{2mm} 

Now take $f\in H^{s+\beta}(\Omega_e)$; by the above observations, $g:=\gamma^{1/2}f \in H^{s+\beta}(\Omega_e)$, and so there exists a unique $w_g\in H^{s+\beta}$ satisfying $((-\Delta)^s+q_\gamma)w=0$ in $\Omega$, $w|_{\Omega_e}=g$. This was proved in \cite{GSU2017}, Lemma 3.1, making use of earliest results found in \cite{VE65}, \cite{Gr14} and \cite{Gr15}. Now let $u_f := \gamma^{-1/2}w_g$. Again by the above observations we have $u_f\in H^{s+\beta}(\Omega_e)$, and by Theorem 3.1 $u_f$ is the unique solution of $\mathrm C^s_\gamma u=0, u|_{\Omega_e}=f$. We also have
\begin{equation*}
\begin{split}
\|\mathrm C^s_\gamma u_f\|_{H^{\beta-s}} & = \| \gamma^{1/2}((-\Delta)^s+q_\gamma)w_g \|_{H^{\beta-s}} \\ & \leq \| \gamma^{1/2} (-\Delta)^s w_g \|_{H^{\beta-s}} + \| w_g (-\Delta)^sm  \|_{H^{\beta-s}} < \infty\;,
\end{split}
\end{equation*}
\noindent and moreover, if $e_h \in H^{s+\beta}(\mathbb R^n)$ is any extension of a given $h\in H^{s+\beta}(\Omega_e)$, 
\begin{equation*}
\begin{split}
\langle \Lambda^s_\gamma f,h \rangle & = B^s_\gamma[u_f, e_h] = \int_{\mathbb R^n}\int_{\mathbb R^n} \nabla^s e_h \cdot (\Theta \cdot \nabla^s u_f) \,dydx = \langle \mathrm C^s_\gamma u_f, e_h \rangle\;.
\end{split}
\end{equation*}
\noindent Given an open set $U$ and a function $u$, let $r_U u := u|_U$. The statement would be proved if we could decompose
$$ \langle \mathrm C^s_\gamma u_f, e_h \rangle = \langle \mathrm r_\Omega C^s_\gamma u_f, r_\Omega e_h \rangle_\Omega + \langle r_{\Omega_e}\mathrm C^s_\gamma u_f, r_{\Omega_e}e_h \rangle_{\Omega_e} \;,$$
\noindent because then since $u_f$ solves the fractional conductivity equation in $\Omega$ we would be able to conclude $\langle \Lambda^s_\gamma f,h \rangle = \langle r_{\Omega_e}\mathrm C^s_\gamma u_f, h \rangle_{\Omega_e}$. In order to use the above decomposition we need to find an $\alpha\in(-1/2, 1/2)$ such that $C^s_\gamma u_f \in H^\alpha$ and $e_h \in H^{-\alpha}$, as in the proof of Lemma 3.1 in \cite{GSU2017}; this task is easily accomplished by taking $\alpha = \beta - s$.  \qed

\vspace{3mm}
\paragraph{\bf Two inverse problems} The two main {\color{black} uniqueness} results about the Calder\'on problem for the fractional Schr\"odinger equation are \cite{RS2017}, Th. 1.1, and \cite{GRSU2018}, Th. 1:  

\paragraph{\bf \emph{Injectivity (infinitely many measurements)}} \emph{Let $\Omega\subset \mathbb R^n$, $n\geq 1$, be bounded open, let $s\in(0,1)$, and let $q_1, q_2 \in L^{n/2s}(\mathbb R^n)$ be such that $0$ is not an eigenvalue of $(-\Delta)^s+q_j$. Let also $W_1, W_2\subset \Omega_e$ be open. If the DN maps for the equations $((-\Delta)^s+q_j)u=0$ in $\Omega$ satisfy} $$ \Lambda_{q_1} [f]|_{W_2} =\Lambda_{q_2} [f]|_{W_2}, \;\;\; \forall f \in C^\infty_c(W_1)\;,$$ \noindent\emph{then $q_1 = q_2$ in $\Omega$.}

\paragraph{\bf \emph{Uniqueness and reconstruction (single measurement)}} \emph{Let $\Omega\subset \mathbb R^n$, $n\geq 1$, be bounded open, let $s\in(0,1)$, and suppose that $0$ is not an eigenvalue of $(-\Delta)^s+q$. Let also $W_1, W_2\subset \Omega_e$ be open, with $\overline\Omega\cap\overline {W_1} =\emptyset$. Assume that either \begin{itemize} \item $s\in[\frac{1}{4},1)$ and $q\in L^{\infty}(\Omega)$, or \item $q\in C^0(\overline\Omega)$. \end{itemize} Given any fixed function $g\in \tilde H^s(W_1)\setminus \{0\}$, the potential $q$ is uniquely determined and can be reconstructed from the knowledge of $\Lambda_{q} [g]|_{W_2}$.}

\vspace{3mm}
\noindent By using the results stated above, one can prove Theorems 1.1 and 1.2.

\paragraph{Proof of Theorem 1.1} If $W_1\cap W_2\neq \emptyset$, there still exist two open sets $W_1'\subset W_1$ and $W_2'\subset W_2$ such that $W_1'\cap W_2'= \emptyset$; so without loss of generality assume that $W_1$, $W_2$ and $\Omega$ are three pairwise disjoint open sets.
\vspace{2mm}

Let $v\in C^\infty_c(W_2)$; the hypothesis of the theorem then reads $$\Lambda^s_{{\gamma_1}}[f]([v])=\Lambda^s_{{\gamma_2}}[f]([v]), \;\;\;\;\; \mbox{for } f\in C^\infty_c(W_1)\;.$$
\noindent Since $\gamma_1=\gamma_2=1$ in $\Omega_e$, one has $\gamma_1^{-1/2}f=\gamma_2^{-1/2}f=f$ in all of $\mathbb R^n$. Therefore, from the previous equality and from Lemma 3.4
\begin{equation*}\begin{split} \Lambda_{q_{\gamma_1}}[f]([v]) & = \Lambda^s_{\gamma_1} [f]|([v]) + \int_{\Omega_e}fv(-\Delta)^sm_1 \, dx \\ & = \Lambda^s_{\gamma_1} [f]([v]) =\Lambda^s_{\gamma_2} [f]([v]) = \Lambda_{q_{\gamma_2}}[f]([v])\;, \end{split}\end{equation*}
\noindent where the integral disappears because supp$(f) \cap$ supp$(v)=\emptyset$. Hence 
\begin{equation}\label{eq:3.6.1} \Lambda_{q_{\gamma_1}}[f]|_{W_2} = \Lambda_{q_{\gamma_2}}[f]|_{W_2},\;\;\;\;\; \mbox{for } f\in C^\infty_c(W_1)\;.\end{equation}
\noindent It is known that $(-\Delta)^sm_j \in L^{n/2s}(\mathbb R^n)$. Therefore, $$ \|q_{\gamma_j}\|^{n/2s}_{L^{n/2s}(\mathbb R^n)} = \int_{\mathbb R^n} \left|\frac{(-\Delta)^sm_j}{\gamma_j^{1/2}}  \right|^{n/2s}dx \leq \underline{\gamma}_j^{-n/4s}\|(-\Delta)^sm_j\|^{n/2s}_{L^{n/2s}(\mathbb R^n)} < \infty\;. $$ Using this and condition (\ref{eq:3.6.1}), one gets $q_{\gamma_1} = q_{\gamma_2}$ in $\Omega$ by the previously stated injectivity result with infinitely many measurements.
\vspace{2mm}

Now let $\bar m = m_2-m_1$; of course supp$(\bar m)\subset \Omega$, and in $\Omega$
\begin{equation}
\label{eq:3.6.2}
\begin{split}
0& = \gamma_1^{1/2}\gamma_2^{1/2}(q_{\gamma_1} - q_{\gamma_2}) = \gamma_1^{1/2}(-\Delta)^sm_2 - \gamma_2^{1/2}(-\Delta)^sm_1 \\ & = (-\Delta)^sm_2 - (-\Delta)^sm_1+m_1(-\Delta)^sm_2 - m_2(-\Delta)^sm_1 \\ & = (1+m_1)(-\Delta)^s\bar m - \bar m (-\Delta)^sm_1 \;.
\end{split}
\end{equation}
\noindent Formula (\ref{eq:3.6.2}) can be written as $(-\Delta)^s\bar m -\frac{(-\Delta)^sm_1}{1+m_1} \bar m =0 $, which shows that $\bar m$ solves the following Dirichlet problem for the fractional Schr\"odinger equation:
$$ \left\{\begin{array}{lr}
        (-\Delta)^s u -\frac{(-\Delta)^sm_1}{1+m_1} u=0 & \text{in } \Omega\,\;\\
        u=0 & \text{in } \Omega_e
        \end{array}\right. \;.$$

\noindent Observe that the equation that $u$ must satisfy in $\Omega$ is the fractional conductivity equation with conductivity $\gamma_1$, by Theorem 3.1. Thus the problem above is well-posed, and so $\bar m =0$ in $\Omega$. This in turn implies $m_1=m_2$, which is the same as saying $\gamma_1 = \gamma_2$ in $\Omega$.
\hfill \qed
\vspace{3mm} 

\paragraph{Proof of Theorem 1.2} By reasoning as before, $W_1$ and $W_2$ can be again supposed to be disjoint. If $v\in H^s(W_2)$, by Lemma 3.4 $$ \Lambda_{q_\gamma} [f]([v]) = \int_{\Omega_e} fv(-\Delta)^sm \, dx +\Lambda_\gamma^s[f]([v]), \;\;\; \forall f\in H^s(\mathbb R^n)\;,$$ \noindent so that, by taking $f=\gamma^{1/2}g$, $$\Lambda_{q_\gamma} [\gamma^{1/2}g]([v]) = \Lambda_\gamma^s[g]|_{W_2} ([v])\;.$$ \noindent Hence $\Lambda_{q_\gamma} [\gamma^{1/2}g]|_{W_2}$ is completely known from $\Lambda_\gamma^s[g]|_{W_2}$. Fix $\epsilon >0$ and observe that the condition $m\in W^{2s+\epsilon,p}_c(\Omega), \forall p> n/\epsilon$ implies $m\in W^{2s, n/2s}_c(\Omega)$ and $(-\Delta)^sm \in C^0(\mathbb R^n)$ by Sobolev embedding theorem. Therefore $q_\gamma \in C^0(\overline\Omega)$, and by the previously stated result concerning uniqueness and reconstruction with a single measurement, $q_\gamma$ can be reconstructed uniquely. By the definition of $q_\gamma$, $m$ solves 
$$ \left\{\begin{array}{lr}
        (-\Delta)^s m - q_\gamma m = - q_\gamma & \text{in } \Omega\,\;\\
        m=0 & \text{in } \Omega_e
        \end{array}\right. \;,$$
\noindent and thus $m$ can be recovered by solving the above problem for Schr\"odinger's equation. 
\hfill \qed
\vspace{5mm} 

\section{A limit case} \label{lim}

Now the previous considerations will be extended to the case $s\rightarrow 1^-$. Since for the fractional Laplacian one has $\lim_{s\rightarrow 1^-}(-\Delta)^s u= -\Delta u$ (\cite{hitch}), it is logical to expect something similar for the other non-local operators. The following holds:

\paragraph{\bf Lemma 4.1} \label{4.1} \emph{Let $u\in H^1(\mathbb R^n)$. Then $\lim_{s\rightarrow 1^-}\|\nabla^s u\|_{L^2(\mathbb R^{2n})} = \|\nabla u\|_{L^2(\mathbb R^{n})}$.}

{\color{black} \paragraph{Remark} This result is a special case of the one given in (\cite{BBM01}), namely when $p=2$. However, since our proof is much easier than the one of the general case, we will still include it for completeness. } 

\paragraph{Proof} Given $i \in \mathbb N$, let $u_i \in C^\infty_c(\mathbb R^n)$ be such that $\|u-u_i\|_{H^1(\mathbb R^n)}\leq 1/i$. By the definition of fractional divergence and Lemma 2.1,
\begin{equation}
\label{eq:4.1.1}
\begin{split}
\lim_{s\rightarrow 1^-} \|\nabla^s u\|&^2_{L^2(\mathbb R^{2n})}  = \lim_{s\rightarrow 1^-} \int_{\mathbb R^n} u(-\Delta)^s u \, dx \\ & = \lim_{i\rightarrow\infty}\lim_{s\rightarrow 1^-} \Big( \int_{\mathbb R^n} u(-\Delta)^s (u-u_i) \, dx + \int_{\mathbb R^n} u(-\Delta)^s u_i \, dx \Big)\;.
\end{split}
\end{equation}
\noindent Since the following estimates hold (\cite{GSU2017}),
\begin{equation}
\label{eq:4.1.new}
\begin{split}
\Big|\int_{\mathbb R^n} u(-\Delta)^s (u-&u_i) \, dx\Big|  = \Big|\int_{\mathbb R^n} (-\Delta)^{s/2}u(-\Delta)^{s/2} (u-u_i) \, dx\Big| \\ & \leq \int_{\mathbb R^n} |(-\Delta)^{s/2}u|\,|(-\Delta)^{s/2} (u-u_i)| \, dx \\ & \leq  \|(-\Delta)^{s/2}u\|_{L^2}\, \|(-\Delta)^{s/2}(u-u_i)\|_{L^2} \\ & \leq \|u\|_{H^s}\, \|u-u_i\|_{H^s} \leq \|u\|_{H^1}\, \|u-u_i\|_{H^1} \leq c/i\;,
\end{split}
\end{equation}
\noindent one gets that $\int_{\mathbb R^n} u(-\Delta)^s (u-u_i) \, dx \rightarrow 0$ upon taking the limits. Moreover $(-\Delta)^s u_i \in \bigcap_{k\in\mathbb N} H^k(\mathbb R^n)\subset L^2(\mathbb R^n)$, and so the second integral in (\ref{eq:4.1.1}) is finite by H\"older. Hence 
\begin{equation}\label{eq:4.1.new2}\begin{split}
\lim_{i\rightarrow\infty}\lim_{s\rightarrow 1^-}  \int_{\mathbb R^n} &u(-\Delta)^s u_i \, dx  = \lim_{i\rightarrow \infty} \int_{\mathbb R^n} u \lim_{s\rightarrow 1^-} (-\Delta)^s u_i \, dx \\ & =  -\lim_{i\rightarrow \infty} \int_{\mathbb R^n} u \Delta u_i \, dx = \lim_{i\rightarrow \infty} \int_{\mathbb R^n} \nabla u \nabla u_i \, dx \\ & = \int_{\mathbb R^n} |\nabla u|^2 \, dx + \lim_{i\rightarrow \infty} \int_{\mathbb R^n} \nabla u \nabla (u_i-u) \, dx = \|\nabla u\|^2_{L^2(\mathbb R^{n})}\; ,
\end{split}\end{equation}
\noindent since the last limit is easily shown to equal $0$ by means of H\"older's inequality. The result is obtained by combining (\ref{eq:4.1.1}), (\ref{eq:4.1.new}) and (\ref{eq:4.1.new2}).
\hfill \qed

{\color{black}\paragraph{Remark} It is not always true that $\nabla^s u(x,y) \rightarrow \nabla u(x)\delta(x-y)$ in distributional sense; quite counter-intuitively, $\lim_{s\rightarrow 1^-} \nabla^s u = 0$ in distributional sense for all $u\in C^\infty_c(\mathbb R^n)$. In fact, if $u \in C^\infty_c(\mathbb R^n)$ and $\phi\in C^\infty_c(\mathbb R^{2n})$, then for some $n$-dimensional balls $B_1, B_2, B_3$ centered at the origin, 

\begin{equation*}
\begin{split}
|\langle \nabla^s u, \phi \rangle| & \leq \int_{\mathbb R^{2n}} |\phi(x,y)| \,|\nabla^su(x,y)|\, dxdy =  \int_{\mathbb R^{2n}} |\phi(x,y)| \,\frac{\mathcal C_{n,s}^{1/2}}{\sqrt{2}}\frac{|u(y)-u(x)|}{|y-x|^{n/2+s}}\, dxdy \\ & \leq c \,\mathcal C_{n,s}^{1/2} \int_{B_1}\int_{B_2} \frac{|u(y)-u(x)|}{|y-x|^{n/2+s}} \, dxdy \leq  c\, \mathcal C_{n,s}^{1/2} \int_{B_1}\int_{B_2} \frac{1}{|y-x|^{n/2+s-1}} \, dxdy \\ & \leq c\, \mathcal C_{n,s}^{1/2} \int_{B_1}\int_{B_3} \frac{1}{|z|^{n/2+s-1}} \, dzdy \leq c \, \mathcal C_{n,s}^{1/2}\;.
\end{split}
\end{equation*}

\noindent Since $\mathcal C_{n,s}^{1/2}$ is bounded by a constant which is independent of $s$ and also $\lim_{s\rightarrow 1^-} \mathcal C_{n,s}^{1/2} = 0$, by dominated convergence the computation above implies $$ \langle \lim_{s\rightarrow 1^-} \nabla^s u, \phi \rangle = \lim_{s\rightarrow 1^-} \langle \nabla^s u, \phi\rangle=0\;.$$

\noindent Observe that this computation is valid also for a more general definition of the fractional gradient, namely one in which $\alpha$ is naturally chosen in such a way that \eqref{modgrad} still holds.}
\vspace{2mm}

\noindent Next, some limit results for the non-local conductivity operator and its DN map. In the rest of this section, the function $m$ will be taken from $W^{2, n/2s}_c(\Omega)$, which embeds into the usual $W^{2s,n/2s}_c(\Omega)$.

\paragraph{\bf Lemma 4.2} \label{4.2} \emph{If $u\in H^2(\mathbb R^n)$, $\lim_{s\rightarrow 1^-} \mathrm C_\gamma^s u = \nabla\cdot(\gamma\nabla u)$ in distributional sense.}

\paragraph{Proof} Let $\phi \in C^\infty_c(\mathbb R^n)$. By reducing the conductivity operator to Schr\"odinger's, one is able to write \begin{equation}\label{eq:4.2.1}\begin{split} \lim_{s\rightarrow 1^-} \int_{\mathbb R^n}& \phi(x) (\nabla\cdot)^s(\Theta\cdot\nabla^su)(x)\, dx = \lim_{s\rightarrow 1^-} \int_{\mathbb R^n} \phi\,\mathrm C_\gamma^su \, dx \\ &= \lim_{s\rightarrow 1^-} \int_{\mathbb R^n} \left(\phi\gamma^{1/2}(-\Delta)^sw -\phi\gamma^{1/2}u(-\Delta)^sm\right) \, dx\;. \end{split}\end{equation}

\noindent Observe now that since $\phi\in C^\infty_c(\mathbb R^n)$ and $u\in H^2(\mathbb R^n)$, we have $\phi u \in H^2(\mathbb R^n)$ as well. Moreover, since $s<1$, we certainly have $m\in W^{2, n/2s}_c(\Omega) \cap L^\infty(\mathbb R^n) \subset W^{2, n/2}_c(\Omega) \cap L^\infty(\mathbb R^n)$; this means that $\gamma^{1/2}$ is a Fourier multiplier on $H^2(\mathbb R^n)$, and therefore $w, \gamma^{1/2}u\phi$ and $\gamma^{1/2}\phi$ all belong to $H^2(\mathbb R^n)$. We can compute
\begin{equation}
\label{eq:mw}
\begin{split}
\| (-\Delta)^sm \|_{H^{-2}} & = \left\| \mathcal F^{-1}\left( \frac{|\xi|^{2s}}{1+|\xi|^2}\hat m(\xi) \right)\right\|_{L^2} \leq c \| \mathcal F^{-1} \hat m(\xi) \|_{L^{2}} = c \| m \|_{L^{2}}\;.
\end{split}
\end{equation}
\noindent In fact, it is easily seen that the function $h_s(x) := \frac{x^{2s}}{1+x^2}$ takes {\color{black} values} in $[0,1)$ for all non-negative $x$ and for all $s\in(0,1)$, which makes $h_s$ a Fourier multiplier on $L^2$. Since $m$ belongs to $L^\infty(\mathbb R^n)$ and has compact support, we see that $\| (-\Delta)^sm \|_{H^{-2}} \leq c \| m \|_{L^{2}} < \infty$, i.e. $(-\Delta)^sm \in H^{-2}(\mathbb R^n)$. By using again (\ref{eq:mw}) with $m$ replaced by $w$, we get $\| (-\Delta)^sw \|_{H^{-2}} \leq c \| w \|_{L^{2}}$; since $w\in H^2(\mathbb R^n)$, this leads to $(-\Delta)^sw \in H^{-2}(\mathbb R^n)$. 

\noindent The above discussion lets us rewrite equation (\ref{eq:4.2.1}) in the form
\begin{equation}\label{eq:divisa} \lim_{s\rightarrow 1^-} \langle \phi, (\nabla\cdot)^s(\Theta\cdot\nabla^su) \rangle = \lim_{s\rightarrow 1^-} \langle \phi\gamma^{1/2},(-\Delta)^sw \rangle-\lim_{s\rightarrow 1^-} \langle \phi\gamma^{1/2}u, (-\Delta)^sm \rangle\;.\end{equation}  

\noindent Trivially, $|h_1(x) - h_s(x)| \leq 2$ for all non-negative $x$ and for all $s\in(0,1)$. With this in mind we can compute
\begin{equation*}
\begin{split}
\| (-\Delta)m - (-\Delta)^sm \|_{H^{-2}} & = \left\| \mathcal F^{-1}\left( \frac{|\xi|^{2}-|\xi|^{2s}}{1+|\xi|^2}\hat m(\xi) \right)\right\|_{L^2} \\ & \leq c \| \mathcal F^{-1} \hat m(\xi) \|_{L^{2}} = c \| m \|_{L^{2}} < \infty\;,
\end{split}
\end{equation*}
\noindent which means that 
\begin{equation*}
\begin{split}
\lim_{s\rightarrow 1^{-}} \| -\Delta m - (-\Delta)^sm \|_{H^{-2}} & = \lim_{s\rightarrow 1^{-}} \left\| \mathcal F^{-1}\left( (h_1(x) - h_s(x)) \hat m(\xi) \right)\right\|_{L^2} \\ & = \left\| \mathcal \lim_{s\rightarrow 1^{-}}(h_1(x) - h_s(x)) \hat m(\xi)\right\|_{L^2} =0\;.
\end{split}
\end{equation*}
\noindent Thus $(-\Delta)^sm \rightarrow -\Delta m$ in $H^{-2}(\mathbb R^n)$ as $s\rightarrow 1^-$, and the same proof can be used to show the analogous result for $(-\Delta)^s w$ as well. We can now deduce from equation (\ref{eq:divisa}) that 
$$ \lim_{s\rightarrow 1^-} \langle \phi, (\nabla\cdot)^s(\Theta\cdot\nabla^su) \rangle = \langle \phi\gamma^{1/2},-\Delta w \rangle- \langle \phi\gamma^{1/2}u, -\Delta m \rangle\;. $$

\noindent Performing some elementary vector calculus computation on this last formula the desired result is immediately obtained: \begin{equation*}\lim_{s\rightarrow 1^-} \int_{\mathbb R^n} \phi\,\mathrm C_\gamma^su \, dx = \int_{\mathbb R^n} \phi \nabla\cdot(\gamma\nabla u)\; dx\,. \end{equation*} \qed

\paragraph{\bf Lemma 4.3} \emph{Let $u,v\in H^1(\mathbb R^n)$. Then $\lim_{s\rightarrow 1} B_\gamma^s[u,v]= \int_{\mathbb R^n} \gamma\nabla u\cdot\nabla v\,dx$.}
\vspace{2mm}

\paragraph{Proof} For all $i\in \mathbb N$, let $u_i, v_i \in C^\infty_c(\mathbb R^n)$ be such that $ \|u-u_i\|_{H^1(\mathbb R^n)}\leq 1/i$ and $\|v-v_i\|_{H^1(\mathbb R^n)}\leq 1/i$. Then we can compute
\begin{equation}
\label{eq:Bsgamma}
\begin{split}
\lim_{s\rightarrow 1^-} B^s_\gamma[u,v] = \lim_{i\rightarrow \infty}\lim_{s\rightarrow 1^-}\left( B^s_\gamma[u-u_i,v] + B^s_\gamma[u_i,v-v_i] + B^s_\gamma[u_i,v_i] \right)\;.
\end{split}
\end{equation}

\noindent By H\"older's inequality we see that
\begin{equation*}
\begin{split}
|B^s_\gamma[u-u_i,v]| & = |\langle \nabla^s(u-u_i), \Theta\cdot\nabla^s v \rangle| \leq \| \nabla^s(u-u_i) \|_{L^2} \|\Theta\cdot \nabla^s v \|_{L^2} \\ & \leq \overline\gamma\|u-u_i\|_{H^s}\|v\|_{H^s} \leq \overline\gamma  \|u-u_i\|_{H^1}\|v\|_{H^1}\;,
\end{split}
\end{equation*}

\noindent so that the first term on the right hand side of (\ref{eq:Bsgamma}) vanishes upon taking the limits. The second term behaves similarly, and so we are left with $\lim_{s\rightarrow 1^-} B^s_\gamma[u,v] = \lim_{i\rightarrow \infty}\lim_{s\rightarrow 1^-}B^s_\gamma[u_i,v_i]$. Now apply Lemma 4.2 to deduce that
\begin{equation*}
\begin{split}
\lim_{s\rightarrow 1^-} B^s_\gamma[u,v] & = \lim_{i\rightarrow \infty}\lim_{s\rightarrow 1^-}B^s_\gamma[u_i,v_i] = \lim_{i\rightarrow \infty}\lim_{s\rightarrow 1^-} \langle \nabla^s u_i, \Theta\cdot\nabla^s v_i \rangle \\ & = \lim_{i\rightarrow \infty}\lim_{s\rightarrow 1^-} \langle u_i, \mathrm C^s_\gamma v_i \rangle = \lim_{i\rightarrow \infty} \langle u_i, \nabla \cdot (\gamma\nabla v_i) \rangle \\ & = \lim_{i\rightarrow \infty} \langle \nabla u_i, \gamma\nabla v_i \rangle \;.
\end{split}
\end{equation*}
 
\noindent The result is now recovered by decomposing this term as in (\ref{eq:Bsgamma}) and then applying again H\"older's inequality.
\hfill \qed
\vspace{5mm}

\paragraph{\bf Corollary 4.4} \emph{Let $f,g\in H^1(\mathbb R^n)$. Then $\lim_{s\rightarrow 1^-} \langle\Lambda_\gamma^s[f],[g]\rangle = \int_{\mathbb R^n} \gamma\nabla u_f\cdot\nabla g\,dx$.}

\paragraph{Proof} The result immediately follows from the previous Lemma and from the definition $\langle\Lambda_\gamma^s[f],[g]\rangle = B^s_\gamma [u_f,g]$.
\hfill \qed
\vspace{5mm}

\section{A simple model: the random walk} \label{model}

This section shows how the non-local conductivity equation naturally arises from weighted long jump random walks. This is an extension of \cite{valdi}, where the fractional Laplacian is related to unweighted long jump random walks. 
\vspace{5mm}

Let $h>0,\; \tau=h^{2s},  \; k\in \mathbb{Z}^n$, $x\in h\mathbb{Z}^n$ and $t\in \tau\mathbb{Z}$. Consider a random walk on the lattice $h\mathbb{Z}^n$, subject to discrete time steps belonging to $\tau\mathbb{Z}$. Define
 \begin{equation*}
  f(x,k) :=
  \begin{cases}
    \gamma^{1/2}(x+hk) |k|^{-n-2s}       & \mbox{if} \quad k\neq 0 \\
    0  & \mbox{if} \quad k=0
  \end{cases}\;. 
	\end{equation*}

\noindent Observe that, $\forall x\in h\mathbb{Z}^n$,
\begin{equation}
\label{eq:sumbound}
\begin{split}
\sum_{k\in\mathbb{Z}^n} f(x,k) & = \sum_{k\in\mathbb{Z}^n \setminus \{0\}} f(x,k) = \sum_{k\in\mathbb{Z}^n \setminus \{0\}} \gamma^{1/2}(x+hk) |k|^{-n-2s} \\ & \leq \overline\gamma^{1/2} \sum_{k\in\mathbb{Z}^n \setminus \{0\}} |k|^{-n-2s} < \infty,
\end{split}
\end{equation}
\noindent and therefore it makes sense to define a normalized version of $f(x,k)$, namely
\begin{equation*}
  P(x,k) :=
  \begin{cases}
    \left( \sum_{j\in\mathbb{Z}^n} f(x,j) \right)^{-1} \gamma^{1/2}(x+hk) |k|^{-n-2s}       & \mbox{if} \quad k\neq 0 \\
    0  & \mbox{if} \quad k=0
  \end{cases}\;. 
	\end{equation*}
\noindent Of course one has $0\leq P(x,k)\leq 1$, and from the definition it follows that 
\begin{equation}
\label{eq:normp}
\begin{split}
\sum_{k\in\mathbb{Z}^n} P(x,k) & = \sum_{k\in\mathbb{Z}^n \setminus \{0\}} P(x,k) = \frac{\sum_{k\in\mathbb{Z}^n \setminus \{0\}} \gamma^{1/2}(x+hk) |k|^{-n-2s}}{\sum_{j\in\mathbb{Z}^n \setminus \{0\}} \gamma^{1/2}(x+hj) |j|^{-n-2s}} =1 \,. 
\end{split}
\end{equation}
\vspace{2mm}

$P(x,k)$ is the probability that a particle found at point $x+hk$ will jump to $x$ in the next discrete step. With $\gamma\equiv 1$ one recovers the case \cite{valdi}, where the probability only depends on the distance between the two points. A non constant function $\gamma$ can instead account for spatially changing properties of the medium, so that the jumping probability is higher from a point whose conductivity is large, while still decreasing with distance. 
\vspace{2mm}

Let $u(x,t)$ be the probability that at some instant $t$ the particle is found at point $x$. It is clearly related to the previous state of the particle by the equation \begin{equation*} u(x, t+\tau) = \sum_{k\in\mathbb{Z}^n\setminus\{0\}}P(x,k)u(x+hk,t) \;\;.  \end{equation*}
\noindent Now compute the time derivative of $u(x,t)$:
\begin{equation*}
\begin{split}
\partial_t u(x,t) & = \lim_{\tau\rightarrow 0} \frac{u(x,t+\tau)-u(x,t)}{\tau} \\ & = \lim_{h\rightarrow 0}\frac{1}{h^{2s}} \left( \sum_{k\in\mathbb{Z}^n\setminus\{0\}}P(x,k)u(x+hk,t) - u(x,t) \right) \\ & = \lim_{h\rightarrow 0}\frac{1}{h^{2s}} \sum_{k\in\mathbb{Z}^n\setminus\{0\}}P(x,k) \left( u(x+hk,t) - u(x,t)\right) \;,
\end{split}
\end{equation*}
 
\noindent where the last line is due to the normalization property (\ref{eq:normp}) of $P(x,k)$. So, \begin{equation}
\label{eq:longguy}
\begin{split}
\partial_t u(x,t) & = \lim_{h\rightarrow 0}\frac{\sum_{k\in\mathbb{Z}^n\setminus\{0\}}\left[\gamma^{1/2}(x+hk) |k|^{-n-2s} \left( u(x+hk,t) - u(x,t)\right)\right]}{h^{2s} \; \sum_{j\in\mathbb{Z}^n \setminus \{0\}} \gamma^{1/2}(x+hj) |j|^{-n-2s}}\;.
\end{split}
\end{equation}

\noindent The denominator is finite, as observed in (\ref{eq:sumbound}), and also bounded away from 0:
\begin{equation}
\label{eq:sumdenom}
\sum_{k\in\mathbb{Z}^n \setminus \{0\}} \gamma^{1/2}(x+hk) |k|^{-n-2s} \geq \underline\gamma^{1/2} \sum_{k\in\mathbb{Z}^n \setminus \{0\}} |k|^{-n-2s} >0\;.
\end{equation}

\noindent By using (\ref{eq:sumdenom}) in equation (\ref{eq:longguy}), one can compute 
\begin{equation*}
\begin{split}
\partial_t u(x,t) & = \lim_{h\rightarrow 0}\frac{\sum_{k\in\mathbb{Z}^n\setminus\{0\}}\left[ h^n \gamma^{1/2}(x+hk) |hk|^{-n-2s} \left( u(x+hk,t) - u(x,t)\right)\right]}{\sum_{j\in\mathbb{Z}^n \setminus \{0\}} \gamma^{1/2}(x+hj) |j|^{-n-2s}} \\ & = C\int_{\mathbb{R}^n} \frac{\gamma^{1/2}(x+z)}{|z|^{n+2s}}\left( u(x+z,t)-u(x,t) \right) dz \\ & = \frac{C}{\gamma(x)^{1/2}} \int_{\mathbb{R}^n} \frac{\gamma^{1/2}(x)\gamma^{1/2}(y)}{|y-x|^{n+2s}}\left( u(y,t)-u(x,t) \right) dy \;,
\end{split}
\end{equation*}

\noindent because the sum approximates the Riemannian integral. Eventually, $\partial_t u(x,t) = \frac{C}{\gamma(x)^{1/2}}\,\mathrm C_\gamma^su$. If $u(x,t)$ is independent of $t$, the fractional conductivity equation $\mathrm C_\gamma^su=0$ is recovered. 
\vspace{5mm}

\section*{References}
\bibliography{Biblio}

\end{document}